\documentclass{article}
\usepackage{amsmath,amssymb,indentfirst}
\usepackage[latin1]{inputenc}
\usepackage{url}
\newtheorem{theorem}{Theorem}
\newtheorem{lemma}{Lemma}

\newtheorem{rmk}{Remark}

\newcommand{\qed}{\hfill\mbox{\raggedright $\Box$}\medskip}

\begin{document}

\title{THE IMPLICIT FUNCTION THEOREM FOR MAPS THAT ARE ONLY DIFFERENTIABLE: AN ELEMENTARY PROOF}

\author{Oswaldo Rio Branco de Oliveira}
\date{}
\maketitle

\begin{abstract} This article shows a very elementary and straightforward proof of the Implicit Function Theorem for differentiable maps $F(x,y)$ defined on a finite-dimensional Euclidean space.
There are no hypotheses on the continuity of the partial derivatives of $F$. The proof employs determinants theory, the mean-value theorem, the intermediate-value theorem, and Darboux's property (the intermediate-value property for derivatives). The proof avoids compactness arguments, fixed-point theorems, and integration theory. A stronger than the classical version of the Inverse Function Theorem is also shown. An example is given.
\end{abstract}

\vspace{0,2 cm}

\hspace{- 0,6 cm}{\sl Mathematics Subject Classification: 26B10, 
26B12}

\hspace{- 0,6 cm}{\sl Key words and phrases:} Implicit Function Theorems, Jacobians, Transformations with Several Variables, Calculus of Vector Functions.

\vspace{0,3 cm}

\section{Introduction.}

The aim of this article is to present a very elementary and straightforward proof of  a version of the Implicit Function Theorem that is  fairly stronger than the classical version. We prove the implicit function theorem for differentiable maps $F(x,y)$, defined on a finite-dimensional Euclidean space, assuming that all the leading principal minors of the partial Jacobian matrix $\frac{\partial F}{\partial y}(x,y)$ are nowhere vanishing
(these hypothesis are already enough to show the existence of an implicit solution)
plus an additional non-degeneracy condition on the matrix $\frac{\partial F}{\partial y}$ to ensure the uniqueness of the implicit solution. 
There are no hypotheses on the continuity of the partial derivatives of the map $F$.

\vspace{0,2 cm}

The results in this article extend de Oliveira \cite{Oliveira1} and \cite{Oliveira2}. In de Oliveira \cite{Oliveira1} are proven the classical versions (enunciated for maps of class $C^1$ on an open set) of the implicit and inverse function theorems. In de Oliveira \cite{Oliveira2} is proven the implicit function theorem for maps $F(x,y)$ such that the partial Jacobian matrix $\frac{\partial F}{\partial y}(x,y)$ is only continuous at the base point.
 
\vspace{0,2 cm}

The proof of the implicit function theorem shown in this article follows Dini's inductive approach (see \cite{Dini}). Moreover, the  proofs of the implicit and the inverse function theorems that we present avoid compactness arguments (i.e., Weierstrass's theorem on maxima), fixed-point theorems (e.g., Banach's fixed point theorem and Brouwer's fixed-point theorem), and Lebesgue's theories of measure and integration. Instead of such tools, we give elementary proofs that are based on the intermediate-value and the mean-value theorems, both on the real line, the intermediate-value property for derivatives on $\mathbb R$ (also known as Darboux's property), and determinants theory.

\vspace{0,2 cm}

 As a corollary of the implicit function theorem shown in this article we obtain a version of the  inverse function theorem that  is stronger than  the classical one proven in most textbooks. An example is given.


\vspace{0,2 cm}
%
%

Some remarks are appropriate regarding proofs of the classical implicit and inverse function theorems. 
Most of these proofs start
 by showing the inverse function theorem  and then derive the implicit function theorem as a trivial consequence. In general, these proofs employ either a compactness argument or the contraction mapping principle (Banach's fixed point theorem), see Krantz and Parks [5, pp. 41--52] and Rudin [6, pp 221--228].
On the other hand, proofs of the classical implicit and inverse function theorems that do not use either of these two tools can be seen in de Oliveira \cite{Oliveira1}.

\vspace{0,2 cm}

Taking into account maps that are everywhere differentiable (their differentials may be everywhere discontinuous), a proof of the implicit function theorem can be found in Hurwicz and Richter \cite{Hurwicz}, whereas a proof of the inverse function theorem can be seen in Saint Raymond \cite{Raymond}. While these two results are quite general, they also have proofs that are quite technical and not that easy to follow. The first of these proofs employs Brouwer's fixed-point theorem  while the second relies on Lebesgue's theories of measure and integration. 

\vspace{0,2 cm}



%
%
%


\vspace{0,1 cm}

Henceforth, we shall freely assume that all the functions are defined on a subset of a finite-dimensional Euclidean space.

\section{Notations and Preliminaries.}

Apart from the intermediate-value and the mean-value theorems, both on the real line, we assume the intermediate-value theorem for derivatives on $\mathbb R$ (also called Darboux's property) stated right below.

\begin{lemma} {\bf (Darboux's Property).} \label{L1} Given $f:[a,b]\to \mathbb R$ differentiable, the image of the derivative function is an interval.
\end{lemma}


Given a $n\times n$ real matrix $A$, we denote its determinant by $\det A$.
The determinant of the sub-matrix of $A$ obtained by deleting the last $n - k$ rows and the last $n-k$ columns of $A$ is the kth order leading principal minor of $A$. 

Given a nonempty subset $X$ of $\mathbb R^n$ and a nonempty subset $Y$ of $\mathbb R^m$, it is well-known that  the cartesian product $X\times Y=\{(x,y): x \in X \ \textrm{and}\ y \in Y\}$
is open in $\mathbb R^n\times \mathbb R^m$ if and only if $X$ and $Y$ are open sets. 

\vspace{0,1 cm}

Let us consider $n$ and $m$, both in $\mathbb N$, and fix the canonical bases $\{e_1,\ldots ,e_n\}$ and $\{f_1, \ldots ,f_m\}$, of $\mathbb R^n$ and $\mathbb R^m$, respectively. Given $x=(x_1,\ldots, x_n)$ and $y=(y_1,\ldots,y_n)$, both in $\mathbb R^n$, we have the inner product $\left<x, y\right>= x_1y_1 + \cdots + x_ny_n$ and 
the norm $|x|=\sqrt{\left<x,x\right>}$. We denote the open ball centered at a point $x$ in $\mathbb R^n$, with radius $r>0$, by 
$B(x;r)=\{y\ \textrm{in}\ \mathbb R^n:\ |y-x|<r\}$. 

\vspace{0,1 cm}

We identify a linear map $T:\mathbb R^n \rightarrow \mathbb R^m$ with the $m\times n$ matrix  $M=(a_{ij})$, where $T(e_j)= a_{1j}f_1 + \cdots + a_{mj}f_m$ for  each $j=1,\ldots,n$. We also write $Tv$ for $T(v)$, where $v\in \mathbb R^n$.


\vspace{0,1 cm}

In this section, $\Omega$ denotes a nonempty open subset of $\mathbb R^n$, where $n\geq 1$. Given a map $F: \Omega \rightarrow \mathbb R^m$ and a point $p$ in $\Omega$, we write $F(p)=\big(F_1(p),\ldots,F_m(p)\big)$. Let us suppose that $F$ is differentiable at $p$. 
The Jacobian matrix of $F$ at $p$ is 
\[JF(p)=\left(\frac{\partial F_i}{\partial x_j}(p)\right)_{\substack{1\leq i\leq m\\1 \leq j\leq n}}=
\left(\begin{array}{lllll}
\frac{\partial F_1}{\partial x_1}(p) & \cdots & \frac{\partial F_1}{\partial x_n}(p)\\
\ \ \ \vdots  &   & \ \ \ \vdots \ \ \      \\
\frac{\partial F_m}{\partial x_1}(p) & \cdots & \frac{\partial F_m}{\partial x_n}(p)
\end{array}
\right).\]
If $F$ is a real function, then we have $JF(p)= \nabla F(p)$, the gradient of $F$ at $p$.

\vspace{0,1 cm}

%

Given $p$ and $q$, both in $\mathbb R^n$, we denote the linear segment with endpoints $p$ and $q$ by $\overline{pq}=\{p+t(q-p):0\leq t\leq 1\}$.
The following result is a trivial corollary of the mean-value theorem on the real line and thus we omit the proof.
\begin{lemma}{\bf (The mean-value theorem in several variables).} \label{L2} Let us consider a differentiable real function $F:\Omega \to \mathbb R$, with $\Omega$ open in $\mathbb R^n$. Let $p$ and $q$ be points in $\Omega$ such that the segment $\overline{pq}$ is within $\Omega$. Then, there exists $c$ in $\overline{pq}$, with $c \neq p$ and $c\neq q$, that satisfies
\[F(p)-F(q) = \left<\nabla F(c), p-q\right>.\]
\end{lemma}

Given a real function $F:\Omega\to \mathbb R$, a short computation shows that the following definition of differentiability is equivalent to that which is most commonly employed. We say that $F$ is differentiable at $p$ in $\Omega$ if there are an open ball $B(p;r)$ within $\Omega$, where $r>0$,  a vector $v$ in $\mathbb R^n$, and a vector-valued  map $E:B(0;r)\to \mathbb R^n$ satisfying
$$ \left\{\begin{array}{ll}
F(p +h)= F(p) + \left<v,h\right> + \left<E(h),h\right>,\ \textrm{for all}\ h \in B(0;r),\\
\textrm{where}\ E(0)=0\ \textrm{and}\  E(h)\to 0\ \textrm{as}\ h\to 0.
\end{array}
\right.
$$


\section{Example and Motivation.}

Right below we give an example of a function $F:\mathbb R^2 \to \mathbb R^2$ so that
$$\left\{\begin{array}{llll}
F\ \textrm{is differentiable everywhere},\\ 
\textrm{the Jacobian matrix $JF$ is not continuous at the origin},\\
\textrm{the leading principal minors of $JF$ do not vanish near the origin},\\
F\ \textrm{is invertible near the origin (proven in the last section)}.
\end{array}
\right.
$$
{\bf Example.} Let us consider the function 
$$F(x,y)=\left\{\begin{array}{ll}
\left(8 x + x^3 \cos \frac{1}{x^2 +y^2},8 y + y^3\sin\frac{1}{x^2 + y^2}\right) & \textrm{outside the origin}\\
(0,0) & \textrm{at the origin}.
\end{array}
\right.
$$
The Jacobian matrix of $F$ outside the origin is given by
$$\left(\begin{array}{lll}
8 + 3x^2 \cos\frac{1}{x^2 +y^2} + \frac{2x^4}{(x^2 +y^2)^2}\sin\frac{1}{x^2 + y^2} & \ \   \frac{2x^3y}{(x^2+y^2)^2}\sin \frac{1}{x^2+y^2}\\
\\
\ \   -\frac{2xy^3}{(x^2+y^2)^2}\cos \frac{1}{x^2+y^2} & 8 + 3y^2 \sin\frac{1}{x^2 +y^2} - \frac{2y^4}{(x^2 +y^2)^2}\cos\frac{1}{x^2 + y^2}\\
\end{array}
\right).$$
On the other hand, a short computation shows that
$$JF(0,0)=\left(\begin{array}{ll}
8 & 0\\
0 & 8\\
\end{array}
\right).
$$

Let us show that $F$ is differentiable at the origin (and thus over the plane). Let $T$ be the linear map associated to the matrix $JF(0,0)$. Given a non-null vector $v=(h,k)$ in the plane we have
$$\frac{F(v) - F(0) - Tv}{|v|}= \frac{\left(8 h + h^3 \cos \frac{1}{h^2 +k^2},8 k + k^3\sin\frac{1}{h^2 + k^2}\right) - (8 h,8 k)}{\sqrt{h^2 + k^2}}\ \ \ \ \ $$
$$\ \ \ \ \ \ \ \ \ \ \ \ \ \ \ \ \ \   =\frac{\left(h^3 \cos \frac{1}{h^2 +k^2}, k^3\sin\frac{1}{h^2 + k^2}\right)}{\sqrt{h^2 + k^2}}\xrightarrow{(h,k)\to (0,0)}(0,0).$$
Thus, $F$ is differentiable at the origin.

We claim that the four entries of $JF(x,y)$ are discontinuous at the origin. For instance, let us look the first entry, which has three terms. The first two terms are continuous at the origin. However, the third term is not. In fact, by employing polar coordinates and writing $(x,y)=(r\cos \theta,r\sin \theta)$ we find
$$\frac{x^4}{(x^2+y^2)^2}\sin\frac{1}{x^2 + y^2} = (\cos^4\theta)\sin \frac{1}{r^2}.$$
Thus, the Jacobian matrix $JF$ is not continuous at the origin.

At last, let us fix $(x,y)$ with $x^2 + y^2\leq 1$.  There exist  $a,b,c,d,e$ and $f$, all in $[-1,1]$, such that the two principal minors of $JF(x,y)$ respectively satisfy
$$|8+ 3a+2b|\geq 3 \ \textrm{and}\ |\det JF(x,y)|=|\left|\begin{array}{ll}
8 + 3a + 2b & 2c\\
2d & 8 + 3e + 2f\\
\end{array}
\right||\geq 3^2 -2^2.
 $$
So, the principal minors of $JF$ do not vanish in the unit disc centered at $(0,0)$.

In the last section we prove that $F$ is invertible on a neighborhood of $(0,0)$.

\newpage

\section{The Implicit  Function Theorem.}


The first implicit function result we prove concerns one equation, several real variables and a differentiable real function.
In its proof, we denote the variable in $\mathbb R^{n+1}=\mathbb R^n\times \mathbb R$ by $(x,y)$, where $x=(x_1,\ldots,x_n)$ is in $\mathbb R^n$ and $y$ is in $\mathbb R$.

\vspace{0,2 cm}

In the next theorem, $\Omega$ denotes a nonempty open set within $\mathbb R^n\times \mathbb R$.

\begin{theorem}\label{TEO1} Let $F:\Omega \to \mathbb R$ be differentiable, with $\frac{\partial F}{\partial y}$ nowhere vanishing, and $(a,b)$ a point in $\Omega$ such that $F(a,b)=0$. Then, there exists an open 
set $X\times Y$, within $\Omega$ and containing the point $(a,b)$, that satisfies the following.
\begin{itemize}
\item[$\bullet$] There exists a unique function $g:X\to Y$ that satisfies $F\big(x,g(x)\big)=0$, for all $x$ in $X$. 
\item[$\bullet$] We have $g(a)=b$. The function $g:X \to Y$ is differentiable and satisfies  
\[\frac{\partial g}{\partial x_j}(x)= - \frac{\frac{\partial F}{\partial x_j}(x,g(x))}{\frac{\partial F}{\partial y}(x,g(x))},\ \textrm{for all}\ x \ \textrm{in}\ X, \ \textrm{where}\ j=1,\ldots,n.\]
\end{itemize}
Moreover, if $\nabla F(x,y)$ is continuous at $(a,b)$ then $\nabla g(x)$ is continuous at $x=a$.
\end{theorem}
{\bf Proof.} By considering the function $F(x+a, \frac{y}{c}+b)$, with $c=\frac{\partial F}{\partial y}(a,b)$, we may assume that $(a,b)=(0,0)$ and $\frac{\partial F}{\partial y}(0,0)=1$. Next, we split the proof into three parts: existence and uniqueness, continuity at the origin, and differentiability.

\begin{itemize}
\item[$\diamond$]{\sf Existence and Uniqueness.} Let us choose a non-degenerate $(n+1)$-dimensional parallelepiped $X\times [-r,r]$, centered at $(0,0)$ and within $\Omega$, whose edges are parallel to the coordinate axes and $X$ is open. Then, the function $\varphi(y)=F(0,y)$, where $y$ runs over $[-r,r]$, is differentiable with $\varphi'$ nowhere vanishing and $\varphi'(0)=1$. Thus, by Darboux's property we have $\varphi'>0$ everywhere and we conclude that $\varphi$ is strictly increasing. Hence, by the continuity of $F$ and shrinking $X$ (if necessary) we may assume that $F$ is strictly negative at the bottom of the parallelepiped and $F$ is strictly positive at the top of the parallelepiped. That is,
$$F\Big|_{X\times \{-r\}}<0\ \ \ \textrm{and}\ \ \ F\Big|_{X\times\{r\}}>0.$$
As a consequence, having fixed an arbitrary $x$ in $X$, the function 
$$\psi(y)=F(x,y),\ \textrm{where}\ y\in [-r,r],$$
satisfies $\psi(-r)<0<\psi(r)$. Hence, by the  mean-value theorem there exists a point $\eta$ in the open interval $Y=(-r,r)$ such that $\psi'(\eta)=\frac{\partial F}{\partial y}(x,\eta)>0$. Therefore, by Darboux's property we have $\psi'(y)>0$ at every $y$ in $Y$. Thus, $\psi$ is strictly increasing and the intermediate-value theorem yields the existence of a unique $y$, we then write $y=g(x)$, in the open interval $Y$ such that $F(x,g(x))=0$.

\item[$\diamond$]{\sf Continuity at  the origin.} Let $\delta$ satisfy $ 0 < \delta < r$. From above, there exists an open set $\mathcal{X}$, contained in $X$ and containing $0$, such that $g(x)$ is in the interval $(-\delta ,\delta)$, for all $x$ in $\mathcal{X}$. Thus, $g$ is continuous at $x=0$.

\item[$\diamond$]{\sf Differentiability.} From the differentiability 
of the real function $F$ at $(0,0)$, and writing $\nabla F(0,0)=(v,1)\in \mathbb R^n\times \mathbb R$ for the gradient of $F$ at $(0,0)$, it follows that there are functions $E_1:\Omega\to \mathbb R^n$ and $E_2:\Omega \to \mathbb R$ satisfying
$$\left\{\begin{array}{ll}
F( h,k)= \left<v, h\right> + k + \left<E_1(h,k),h\right> + E_2(h,k)k,\\
\\
\textrm{where}\ \lim\limits_{(h,k)\to (0,0)}E_j(h,k)= 0 = E_j(0,0), \ \textrm{for}\ j=1,2.

\end{array}
\right.
$$

Hence, substituting 
[we already proved that $g(h)\xrightarrow{h\to 0}g(0)=0$]
$$\left\{\begin{array}{lll}
k=g(h),\\
E_j\big(h,g(h)\big)=\epsilon_j(h),\ 
\textrm{with}\ \lim\limits_{h\to 0}\epsilon_j(h)=\epsilon_j(0)=0\ \textrm{for}\ j=1,2,
\end{array}
\right.
$$
and noticing that we have $F\big(h,g(h)\big)=0$, for all possible $h$, we obtain
$$\left<v,h\right> + g(h) \ + \left<\epsilon_1(h),h\right> +  \epsilon_2(h)g(h)=0.\ \ \ \ \ \ \ \ \ \ \  $$
Thus, 
$$ [1 + \epsilon_2(h)]g(h) = - \left<v,h\right> - \left<\epsilon_1(h),h\right>. $$
If $|h|$ is small enough, then we have $1 + \epsilon_2(h)\neq 0$ and we may write

$$g(h)= \left<-v,h\right> +\left<\epsilon_3(h),h\right>,$$
where 
$$\epsilon_3(h)= \frac{\epsilon_2(h)}{1+\epsilon_2(h)}v -\frac{\epsilon_1(h)}{1+\epsilon_2(h)} \ \textrm{and}\ \lim_{h\to 0}\epsilon_3(h)=0.$$
Therefore, $g$ is differentiable at $0$ and $\nabla g(0)=-v$. 

Now, given any $a'$ in $X$, we put $b'=g(a')$. Then, $g: X \to Y$ solves the problem $F\big(x,h(x)\big)=0$, for all $x$ in $X$, with the condition $h(a')= b'$. From what we have just done it follows that $g$ is differentiable at $a'$. 
\end{itemize}
$\qed$

Next, we prove the implicit function theorem for a finite number of equations. 
Some notations are appropriate. We denote the variable in 
$\mathbb R^n\times \mathbb R^m = \mathbb R^{n+m}$ by $(x,y)$, where $x=(x_1,\ldots,x_n)$ is in $\mathbb R^n$ and $y=(y_1,\ldots,y_m)$ is in $\mathbb R^m$. 
Given $\Omega$ an open subset of $\mathbb R^n\times \mathbb R^m$ and a  
differentiable map
$F:\Omega \to \mathbb R^m$ we write   $F=(F_1,\ldots,F_m)$, with $F_i$ the ith component of $F$ and $i=1,\ldots, m$, and 
\[\frac{\partial F}{\partial y}=\left(\frac{\partial F_i}{\partial y_j}\right)_{\substack{1\leq i\leq m\\1 \leq j\leq m}}=
\left(\begin{array}{lllll}
\frac{\partial F_1}{\partial y_1} & \cdots & \frac{\partial F_1}{\partial y_m}\\
\ \ \ \vdots  &   & \ \ \ \vdots \ \ \      \\
\frac{\partial F_m}{\partial y_1} & \cdots & \frac{\partial F_m}{\partial y_m}
\end{array}
\right).\]
Analogously, we define the matrix $\frac{\partial F}{\partial x}=\big(\frac{\partial F_i}{\partial x_k}\big)$, where $1\leq i \leq m$ and $1\leq k \leq n$.

\begin{theorem}\label{TEO2}{\bf (The Implicit Function Theorem).}  Let $F:\Omega \to \mathbb R^m$ be differentiable, with $\Omega$ a non-degenerate open ball within $\mathbb R^n\times \mathbb R^m$ and centered at $(a,b)$. Let us suppose that $F(a,b)=0$ and that all the leading principal minors of the matrix $\frac{\partial F}{\partial y}$ are nowhere vanishing. 
The following are true.
\begin{itemize}
\item[$\bullet$] There exists an open 
set $X\times Y$, within $\Omega$ and containing $(a,b)$, and a differentiable function $g:X\to Y$ that satisfies 
$$F\big(x,g(x)\big)=0, \ \textrm{for all} \ x \in \ X,\ \ \textrm{and}\ \ g(a)=b. $$
\item[$\bullet$] We have
\[Jg(x)  = - \left[\frac{\partial F}{\partial y }(x,g(x))\right]_{m\times m}^{-1}\left[\frac{\partial F}{\partial x}(x,g(x))\right]_{m\times n},\ \textrm{for all}\ x \ \textrm{in}\ X.\]
\end{itemize}
Let us suppose that we also have 
$\det\big(\frac{\partial F_i}{\partial y_j}(\xi_{ij})\big)_{1\leq i,j\leq m}\neq 0$, for all $\xi_{ij}$ in $\Omega$ and $1\leq i,j\leq m$.
 Then, the following is true.
\begin{itemize}
\item[$\bullet$] If $h:X\to Y$ satisfies $F\big(x,h(x)\big)=0$ for all $x$ in $X$, then we have $h=g$.
\end{itemize}
\end{theorem}
{\bf Proof.} Let us split the proof into three parts: existence and differentiability, differentiation formula, and uniqueness.

\begin{itemize}
\item[$\diamond$]{\sf Existence and differentiability.} We claim that the system
$$\left\{\begin{array}{ll}
F_1(x,y_1,\ldots,y_m) =0,\\
F_2(x,y_1,\ldots,y_m) =0,\\
\ \ \ \ \ \ \  \ \ \ \ \ \ \vdots\\
F_m(x,y_1,\ldots,y_m) =0,\\ 
\end{array}
\right. \  \textrm{with the conditions}\ \ 
\left\{\begin{array}{llll}
y_1(a)=b_1\\
y_2(a)=b_2\\
\ \ \ \   \ \  \vdots\\
y_m(a)=b_m,\\
\end{array}
\right.
$$
has a differentiable solution $g(x)=\big(g_1(x),\ldots, g_m(x)\big)$ on some open 
set $X$ containing $a$ [i.e., we have $F\big(x,g(x)\big)=0$ for all $x$ in $X$ and $g(a)=b$].

Let us employ induction on $m$. The case $m=1$ follows immediately from Theorem \ref{TEO1}.

Assuming that the claim holds for $m-1$, let us examine the case $m$. 

Then,
 given a pair $(x,y)=(x,y_1,\ldots, y_m)$ we introduce the helpful notations $y'=(y_2,\ldots,y_m)$, $y=(y_1,y')$, and $(x,y)=(x,y_1,y')$.


As a first step, we consider the equation 
$$F_1(x,y_1,y')=0, \ \textrm{with the condition} \ y_1(a,b')=b_1,$$
 where $x$ and $y'$ are independent variables and $y_1$ is the dependent one. 
Since $\frac{\partial F_1}{\partial y_1}(x,y_1,y')$ is nowhere vanishing, 
by Theorem \ref{TEO1} it follows that there exists a differentiable function $\varphi(x,y')$ on some open 
set [let us say, $\mathcal{X}\times \mathcal{Y}'$] containing $(a,b')$ that satisfies
\[F_1[x,\varphi(x,y'),y']=0\ \textrm{(on $\mathcal{X}\times \mathcal{Y}'$)} \  \textrm{and the condition}\ \varphi(a,b')=b_1.\]
From Theorem \ref{TEO1} we see that $\varphi(x,y')$ also satisfies the $m-1$ equations
\begin{equation}\label{EQN1} -\frac{\partial \varphi}{\partial y_j}(x,y') =
\frac{\frac{\partial F_1}{\partial y_j}[x,\varphi(x,y'),y']}{\frac{\partial F_1}{\partial y_1}[x,\varphi(x,y'),y']},
\ \textrm{for all}\ j=2,\ldots,m. 
\end{equation}

As a second step, we look at solving the system with $m-1$ equations
 $$\left\{\begin{array}{lll}
F_2[x,\varphi(x,y'),y']=0\\
\ \ \ \ \ \ \ \ \ \ \ \ \vdots\\
F_m[x,\varphi(x,y'),y']=0\\ 
\end{array}
\right. , 
\  \textrm{with the condition}\ \ 
y'(a)=b'.
$$
Here, $x$ is the independent variable while $y'$ is the dependent variable.
Let us define $\mathcal{F}_i(x,y')=F_i[x,\varphi(x,y'),y']$, with $i=2,\ldots,m$, and  write $\mathcal{F}=(\mathcal{F}_2,\ldots,\mathcal{F}_m)$. Evidently, the map $\mathcal{F}$ is differentiable. In order to employ the induction hypothesis, let us show that  all the leading principal minors of the partial Jacobian matrix $\frac{\partial \mathcal{F}}{\partial y'}$ are nowhere vanishing.

Thus, let us consider the leading principal minor (a general one)
$$\left| \begin{array}{llllllll}
\frac{\partial F_2}{\partial y_1}\frac{\partial \varphi}{\partial y_2} + \frac{\partial F_2}{\partial y_2} \,&\, \frac{\partial F_2}{\partial y_1}\frac{\partial \varphi}{\partial y_3} + \frac{\partial F_2}{\partial y_3} \,&\, \cdots \,&\, 
\frac{\partial F_2}{\partial y_1}\frac{\partial \varphi}{\partial y_k} + \frac{\partial F_2}{\partial y_k}\\
\\
\frac{\partial F_3}{\partial y_1}\frac{\partial \varphi}{\partial y_2} + \frac{\partial F_3}{\partial y_2} \,&\, \frac{\partial F_3}{\partial y_1}\frac{\partial \varphi}{\partial y_3} + \frac{\partial F_3}{\partial y_3} \,&\, \cdots \,&\, 
\frac{\partial F_3}{\partial y_1}\frac{\partial \varphi}{\partial y_k} + \frac{\partial F_3}{\partial y_k}\\

\ \ \ \ \ \ \ \ \ \vdots \,&\, \ \ \ \ \ \ \ \ \ \vdots \,&\, \,&\, \ \ \ \ \ \ \ \ \   \vdots\\

\frac{\partial F_k}{\partial y_1}\frac{\partial \varphi}{\partial y_2} + \frac{\partial F_k}{\partial y_2} \,&\,
 \frac{\partial F_k}{\partial y_1}\frac{\partial \varphi}{\partial y_3} + \frac{\partial F_k}{\partial y_3} \,&\, \cdots \,&\, 
\frac{\partial F_k}{\partial y_1}\frac{\partial \varphi}{\partial y_k} + \frac{\partial F_k}{\partial y_k}\\
\end{array}
\right|.
$$

Developing this determinant by the columns and then canceling the everywhere vanishing determinants we arrive at (a sum of $k$ determinants)
$$\det \left(\frac{\partial \mathcal{F}_i}{\partial y_j}\right)_{2\leq i,j\leq k}= \ \left| \begin{array}{llllll}
\frac{\partial F_2}{\partial y_2} &  \frac{\partial F_2}{\partial y_3} & \cdots & 
\frac{\partial F_2}{\partial y_k}\\
\\
\frac{\partial F_3}{\partial y_2} & \frac{\partial F_3}{\partial y_3} & \cdots & 
\frac{\partial F_3}{\partial y_k}\\

\ \ \vdots & \ \ \vdots & & \ \    \vdots\\

\frac{\partial F_k}{\partial y_2} & \frac{\partial F_k}{\partial y_3} & \cdots & 
\frac{\partial F_k}{\partial y_k}\\
\end{array}
\right|\ \ \ \ \ \ \ \ \ \ \ \ \ \ \ \ \ \ \ \ \ \ \ \ \ \ \ \ \ \ \ \ \ \ \ \ \ \ \ \ \ \ \ \ \ \ $$
$$+ \left| \begin{array}{llllll}
\frac{\partial F_2}{\partial y_1}\frac{\partial \varphi}{\partial y_2} & 
 \frac{\partial F_2}{\partial y_3} & \cdots & 
\frac{\partial F_2}{\partial y_k}\\
\\
\frac{\partial F_3}{\partial y_1}\frac{\partial \varphi}{\partial y_2} & 
\frac{\partial F_3}{\partial y_3} & \cdots & 
\frac{\partial F_3}{\partial y_k}\\

\ \ \ \ \, \vdots & \ \  \vdots & & \ \   \vdots\\

\frac{\partial F_k}{\partial y_1}\frac{\partial \varphi}{\partial y_2} & 
\frac{\partial F_k}{\partial y_3} & \cdots & 
\frac{\partial F_k}{\partial y_k}\\
\end{array}
\right| 
+ \left| \begin{array}{llllll}
 \frac{\partial F_2}{\partial y_2} \,&\, 
\frac{\partial F_2}{\partial y_1}\frac{\partial \varphi}{\partial y_3} \,&\, \frac{\partial F_2}{\partial y_4} \,&\, \cdots \,&\, 
\frac{\partial F_2}{\partial y_k}\\
\\
\frac{\partial F_3}{\partial y_2} \,&\, 
\frac{\partial F_3}{\partial y_1}\frac{\partial \varphi}{\partial y_3} \,&\, \frac{\partial F_3}{\partial y_4} \,&\, \cdots \,&\, 
\frac{\partial F_3}{\partial y_k}\\

\ \  \vdots \,&\, \ \ \ \  \vdots \,& \ \ \vdots \,&\, \,&\, \ \   \vdots\\

\frac{\partial F_k}{\partial y_2} \,&\, 
\frac{\partial F_k}{\partial y_1}\frac{\partial \varphi}{\partial y_3} \,&\, \frac{\partial F_k}{\partial y_4} \,&\, \cdots \,&\, 
\frac{\partial F_k}{\partial y_k}\\
\end{array}
\right|$$
$$+\cdots + \left| \begin{array}{llllll}
\frac{\partial F_2}{\partial y_2} \,&\, 
\cdots \,&\, \frac{\partial F_2}{\partial y_{k-1}} \,&\, 
\frac{\partial F_2}{\partial y_1}\frac{\partial \varphi}{\partial y_k}\\
\\
\frac{\partial F_3}{\partial y_2} \,&\, 
 \cdots \,&\,  \frac{\partial F_3}{\partial y_{k-1}} \,&\, 
\frac{\partial F_3}{\partial y_1}\frac{\partial \varphi}{\partial y_k}\\

\ \  \vdots \,&\, \ \  \,&\, \ \ \vdots \,&\, \ \   \ \  \vdots\\

\frac{\partial F_k}{\partial y_2} \,&\, 
\cdots \,&\, \frac{\partial F_k}{\partial y_{k-1}} \,&\, 
\frac{\partial F_k}{\partial y_1}\frac{\partial \varphi}{\partial y_k}\\
\end{array}
\right|.\ \ \ \ \ \ \ \ \ \ \ \ \ \ \ \ \ \ \ \ \ \ \ \ \ \ \ \ \ \ \ \ \ \ \  \ \  $$
Thus, we obtain (keeping track of $\frac{\partial \varphi}{\partial y_j}$ for  $j$ even and also for $j$ odd)
$$\det \left(\frac{\partial \mathcal{F}_i}{\partial y_j}\right)_{2\leq i,j\leq k}=\left| \begin{array}{llllll}
1 &  -\frac{\partial \varphi}{\partial y_2} &-\frac{\partial \varphi}{\partial y_3} & \cdots & 
-\frac{\partial \varphi}{\partial y_k}\\
\\
\frac{\partial F_2}{\partial y_1} &  \ \ \frac{\partial F_2}{\partial y_2} &\ \ \frac{\partial F_2}{\partial y_3} & \cdots &
\ \ \frac{\partial F_2}{\partial y_k}\\
\\
\frac{\partial F_3}{\partial y_1} & \ \ \frac{\partial F_3}{\partial y_2} & \ \ \frac{\partial F_3}{\partial y_3} & \cdots & 
\ \ \frac{\partial F_3}{\partial y_k}\\

\ \ \vdots & \ \ \ \ \vdots & \ \ \ \ \vdots & & \ \    \ \ \vdots\\

\frac{\partial F_k}{\partial y_1} & \ \ \frac{\partial F_k}{\partial y_2} & \ \ \frac{\partial F_k}{\partial y_3} & \cdots & 
\ \ \frac{\partial F_k}{\partial y_k}\\
\end{array}
\right|.\ \ \ \ \ \ \ \ \ \ \ \ \ \ \ \  \ \ \ \ \ \ \  \ \ \ \ \ \ \ \ \ \ \ \ \  \ \ \ \ \ \ \ \ \ $$
The already remarked identity $-\frac{\partial \varphi}{\partial y_j} = \frac{\partial F_1}{\partial y_j}/\frac{\partial F_1}{\partial y_1}$ [see formula (\ref{EQN1})] leads to
$$\det \left(\frac{\partial \mathcal{F}_i}{\partial y_j}\right)_{2\leq i,j\leq k}=\frac{1}{\frac{\partial F_1}{\partial y_1}}\left| \begin{array}{llllll}
\frac{\partial F_1}{\partial y_1} &  \frac{\partial F_1}{\partial y_2} & \frac{\partial F_1}{\partial y_3} & \cdots & 
\frac{\partial F_1}{\partial y_k}\\
\\
\frac{\partial F_2}{\partial y_1} &   \frac{\partial F_2}{\partial y_2} & \frac{\partial F_2}{\partial y_3} & \cdots & 
 \frac{\partial F_2}{\partial y_k}\\
\\
\frac{\partial F_3}{\partial y_1} &  \frac{\partial F_3}{\partial y_2} &  \frac{\partial F_3}{\partial y_3} & \cdots & 
 \frac{\partial F_3}{\partial y_k}\\

\ \ \vdots & \ \ \vdots & \ \ \vdots & & \ \    \vdots\\

\frac{\partial F_k}{\partial y_1} &  \frac{\partial F_k}{\partial y_2} &  \frac{\partial F_k}{\partial y_3} & \cdots & 
 \frac{\partial F_k}{\partial y_k}\\
\end{array}
\right|.$$
Hence, all the leading principal minors of $\frac{\partial \mathcal{F}}{\partial y'}$ are nowhere vanishing.

Thus, by induction hypothesis there exists a differentiable function $\psi$ defined on an open 
set $X$ containing $a$  [with  $\psi(X)$ within $\mathcal{Y}'$] that satisfies 
\[ \left\{\begin{array}{ll}
F_i[x,\varphi\big(x,\psi(x)\big),\psi(x)\big]=0, \ \textrm{for all}\ x \ \textrm{in}\ X, \ \textrm{for all} \ i=2,\ldots,m,\\ 
\textrm{and the condition} \ \psi(a)=b'.
\end{array}
\right.\] 
Clearly, we also have $F_1\big[x,\varphi\big(x,\psi(x)\big),\psi(x)\big]=0$, for all $x$ in $X$. Defining 
$$g(x) = \big(\varphi(x,\psi(x)),\psi(x)\big), \ \textrm{where}\ x \in X,$$ 
we obtain $F[x,g(x)]=0$, for every  $x$ in $X$, with $g$ differentiable on $X$, and also the identity $g(a)=\big(\varphi(a,b'),b'\big)=(b_1,b')=b$.
\item[$\diamond$] {\sf Differentiation formula.} Differentiating $F[x,g(x)]=0$ we find
\[\frac{\partial F_i}{\partial x_k} + \sum_{j=1}^m\frac{\partial F_i}{\partial y_j}\frac{\partial g_j}{\partial x_k}=0,\ \textrm{with}\ 1\leq i\leq m\ \textrm{and}\ 1\leq k\leq n.\]
In matricial form, we write $\frac{\partial F}{\partial x}\big(x,g(x)\big) + \frac{\partial F}{\partial y}\big(x,g(x)\big)Jg(x)=0$. 

\item[$\diamond$] {\sf Uniqueness.} If $h:X\to Y$ 
and $x$ in $X$ satisfy $F(x,h(x))=0$, by Lemma \ref{L2} (the mean-value theorem in several variables) there exist $c_1,\ldots, c_m$, all in the open ball $\Omega$ (a convex set), satisfying
\begin{align*} 0 & = F\big(x,h(x)\big)-F\big(x,g(x)\big)\\
& =
\left[\begin{array}{llll}
\frac{\partial F_1}{\partial y_1}(c_1) &\cdots & \frac{\partial F_1}{\partial y_m}(c_1)\\
\ \ \ \ \vdots &  & \ \ \ \ \vdots\\
\frac{\partial F_m}{\partial y_1}(c_m) &\cdots & \frac{\partial F_m}{\partial y_m}(c_m)\\
\end{array}
\right]
\left[\begin{array}{llll}
h_1(x)-g_1(x)\\
\ \ \ \ \ \ \ \ \ \vdots\\
h_m(x)-g_m(x)\\
\end{array}
\right]. \\
\end{align*}
From the hypothesis we have $\det\big(\frac{\partial F_i}{\partial y_j}(c_i)\big)\neq 0$. Thus, $h(x)=g(x)$.
\end{itemize}
$\qed$

\section{The Inverse Function Theorem.}


\begin{theorem} \label{TEO3}{\bf (The Inverse Function Theorem).} Let $F:\Omega \to \mathbb R^n$ be differentiable, with $\Omega$ a non-degenerate open ball within $\mathbb R^n$ and centered at $x_0$. Let us suppose that all the leading principal minors of $JF(x)$ are nowhere vanishing. We also suppose  $\det\big(\frac{\partial F_i}{\partial x_j}(\xi_{ij})\big)_{1\leq i,j\leq n}\neq 0$ for all $\xi_{ij} \in \Omega$, where $i,j=1,\ldots,n$. Under such conditions, there exist an open set $X$ containing $x_0$, an open set $Y$ containing $y_0=F(x_0)$, and a differentiable  $G:Y\to X$ satisfying
$$F\big(G(y)\big)=y, \ \textrm{for all}\ y \in Y, \ \textrm{and}\ G\big(F(x)\big)=x, \ \textrm{for all}\ x \in X.$$
 In addition, 
\[JG(y)= JF\big(G(y)\big)^{-1}, \ \textrm{for all} \ y \ \textrm{in}\ Y.\]
\end{theorem}
{\bf Proof.} Let us split it into two parts: injectivity of $F$ and existence of $G$.

\begin{itemize}
\item[$\diamond$] {\sf Injectivity of $F$.} Let us suppose $F(p)=F(q)$, with $p$ in $\Omega$ and $q$ in $\Omega$. By employing Lemma \ref{L2} we obtain $c_1,\ldots,c_n$, all in the ball $\Omega$, such that
$$0=F(p)-F(q)= 
\left[\begin{array}{llll}
\frac{\partial F_1}{\partial x_1}(c_1) &\cdots & \frac{\partial F_1}{\partial x_n}(c_1)\\
\ \ \ \ \vdots &  & \ \ \ \ \vdots\\
\frac{\partial F_n}{\partial x_1}(c_n) &\cdots & \frac{\partial F_n}{\partial x_n}(c_n)\\
\end{array}
\right]
\left[\begin{array}{lll}
p_1-q_1\\
\ \ \ \   \vdots\\
p_n-q_n\\
\end{array}
\right].$$
The hypotheses imply $\det\big(\frac{\partial F_i}{\partial x_j}(c_i)\big)\neq 0$. Thus, $p=q$.

\item[$\diamond$] {\sf Existence of $G$.} The map 
$$\Phi(y,x)= F(x)-y, \ \textrm{where}\ (y,x) \in \mathbb R^n\times \Omega,$$
 is differentiable and $\Phi(y_0,x_0)=0$. From the hypotheses it follows that all the leading principal minors of 
$\frac{\partial \Phi}{\partial x}(y,x)=JF(x)$ are nowhere vanishing in $\mathbb R^n\times \Omega$ and 
$$\det\left(\frac{\partial \Phi}{\partial x}(\eta_{ij},\xi_{ij})\right)= \det\left(\frac{\partial F}{\partial x}(\xi_{ij})\right)\neq 0,$$
for all $(\eta_{ij},\xi_{ij})\in \mathbb R^n\times \Omega$. 
 The Implicit Function Theorem guarantees an open set $Y$ containing $y_0$ and a differentiable map $G:Y \to \Omega$ satisfying  
$$F\big(G(y)\big)=y,\ \textrm{for all}\ y \ \textrm{in}\ Y.$$

Thus, $G$ is bijective from $Y$ to $X=G(Y)$ and $F$ is bijective from $X$ to $Y$. We also have $X=F^{-1}(Y)$. Since $F$ is continuous, the set $X$ is open (and contains $x_0$).

Putting $F(x)=\big(F_1(x),\ldots,F_n(x)\big)$ and $G(y)=\big(G_1(y),\ldots,G_n(y)\big)$ and differentiating  $\big(F_1(G(y)),\ldots,F_n(G(y))\big)$ we find
\[ \sum_{k=1}^n\frac{\partial F_i}{\partial x_k}\frac{\partial G_k}{\partial y_j}= \frac{\partial y_i}{\partial y_j}=
\left\{\begin{array}{ll}
1,\ \textrm{if} \ i=j,\\
0, \ \textrm{if}\ i\neq j.
\end{array}
\right.\]
\end{itemize}
 $\qed$

\

\begin{rmk} Is is clear that the function in section 3 (Example and Motivation) satisfy the conditions of the above inverse function theorem and is thus invertible, with differentiable inverse function, on a neighborhood of the origin.
\end{rmk}
\begin{rmk} It is not difficult to see that Theorem \ref{TEO2} implies the implicit function theorem for a differentiable function $F: \Omega\subset \mathbb R^n\times \mathbb R^m\to \mathbb R^m$, with $F(a,b)=0$ and $\Omega$ an open set, whose partial Jacobian matrix $\frac{\partial F}{\partial y}(x,y)$ is continuous at the base point $(a,b)$ and $\det\frac{\partial F}{\partial y}(a,b)\neq 0$. In fact, by a linear change of coordinates in the $y$ variable, we may assume $\frac{\partial F}{\partial y}(a,b)=I$, with $I$ the $m\times m$ identity matrix. Thus, on some open neighborhood of $(a,b)$, we have $\det\big(\frac{\partial F_i}{\partial y_j}(\xi_{ij})\big)_{1\leq i,j\leq k}\neq 0$ for all $\xi_{ij}$ in this neighborhood, where $1\leq i,j\leq k$, for each $k=1,\ldots,m$.
\end{rmk}
\begin{rmk}
Similarly, Theorem \ref{TEO3} implies the inverse function theorem for a differentiable function $F: \Omega\subset \mathbb R^n\to \mathbb R^n$, with $F(x_0)=y_0$ and $\Omega$ an open set in $\mathbb R^n$, whose Jacobian matrix $JF(x)$ is continuous at $x_0$ and $\det JF(x_0)\neq 0$.
\end{rmk}

\paragraph{Acknowledgments.}  The author is gratful to Professor P.  A. Martin for discussions that lead to the example given in this article.



\bigskip

\noindent\textit{Departamento de Matemática,
Universidade de São Paulo\\
Rua do Matão 1010 - CEP 05508-090\\
São Paulo, SP - Brasil\\
oliveira@ime.usp.br}

\bigskip

\end{document}